\documentclass{arxiv}
\usepackage{amssymb}
\usepackage{amsmath}
\usepackage{graphicx}
\usepackage{epstopdf}
\usepackage{color}
\usepackage{bbold}
\usepackage{units}

\newtheorem{theorem}{Theorem}

\newtheorem{remark}{Remark}
\newtheorem{assumption}{Assumption}

\newtheorem{definition}{Definition}

\newcommand{\ones}{\mathbf{1}}
\newcommand{\tr}{\intercal}

\begin{document}

\begin{frontmatter}

\title{Verifying Control Invariance of H-Polytopes is coNP-Complete}

\author[IMT]{Mario Eduardo Villanueva},
\ead{me.villanueva@imtlucca.it}
\author[ST,IMT]{Boris Houska}
\ead{borish@shanghaitech.edu.cn}

\address[IMT]{IMT School for Advanced Studies Lucca, Italy}%
\address[ST]{ShanghaiTech University, China}

\begin{keyword}
Linear Systems, Control Theory, Computational Complexity
\end{keyword}

\begin{abstract}
We prove that verifying whether an H-polytope is control invariant for a discrete-time linear system with polytopic input constraints is \textsc{coNP}-complete. This establishes a fundamental computational complexity barrier for set-theoretic control. As an immediate consequence, the exact computation of maximal control invariant subsets of H-polytopes is \textsc{coNP}-hard. These hardness results hold regardless of whether the input constraint set is given in H- or V-representation.

\end{abstract}

\end{frontmatter}

\section{Introduction}

A set is control invariant if, starting from every state in it, some admissible control input keeps the successor state inside the set. The notion originated in the theory of control invariant subspaces \cite{Basile1969,Wonham1970}. Control invariant sets are the foundation for a wide variety of set-theoretic methods in control, as surveyed in~\cite{Blanchini1999,Blanchini2015,Houska2025}. For example, in model predictive control, they serve as terminal regions to guarantee recursive feasibility and stability \cite{Mayne2000,Rawlings2015}.

Note that this paper makes no attempt to survey the extensive literature on computing control invariant sets. This is because we ask a much more fundamental question: Given the system matrices $A$ and $B$ of a linear discrete-time system, a polytopic input constraint $\mathbb U$, and a candidate polytope $X$, Can we decide in polynomial time whether $X$ is control invariant?

When the candidate $X$ is a V-polytope--that is, given by its vertices--verification is easy. Namely, Gutman and Cwikel \cite{Gutman1986} showed that it suffices to check the control invariance conditions at the vertices of the given candidate polytope, reducing the problem to a set of linear programs~\cite{Khachiyan1979}. For H-polytopes--specified by halfspace inequalities--their approach is, however, not viable. This is because a polytope in halfspace representation may have exponentially many vertices~\cite{Mullen1971}. Therefore, the question of whether verification of control invariance of H-polytopes is tractable has remained wide open.

\subsection*{Contribution}

In this paper, we settle the above open question by proving that verifying control invariance of an H-polytope is \textsc{coNP}-complete, see Theorem~\ref{thm::main}. Consequently, unless $\textsc{P} = \textsc{NP}$, no polynomial-time verification procedure exists for the general problem. 
The hardness result extends immediately to exact computation of maximal control invariant subsets of given H-polytopes, see Corollary~\ref{cor::MCIset}.

\section{Preliminaries}
This section briefly reviews complexity classes, polytope representations, and the concept of control invariance.

\subsection{Complexity Classes}

A decision problem is in \textsc{coNP} if its ``no'' instances have polynomial-time verifiable certificates; equivalently, its complement is in \textsc{NP}~\cite{Arora2009}. To compare the difficulty of problems, we rely on polynomial-time reductions. Intuitively, if problem $L$ can be reduced to problem $K$ in polynomial time, then $L$ is no harder than $K$: a solver for $K$ can be used to solve $L$ after a suitable input transformation.

Let $L$ be our target problem and $K$ a known \textsc{coNP}-complete problem. We define the relevant notions as follows:
\begin{itemize}
    \item $L$ is \textsc{coNP}-hard if every problem in \textsc{coNP} can be reduced to $L$. In practice, it suffices to establish a polynomial-time reduction from $K$ to $L$; hardness for the entire class then follows by transitivity, since every problem in \textsc{coNP} already reduces to $K$.
    \item $L$ is \textsc{coNP}-complete if it is both \textsc{coNP}-hard and belongs to \textsc{coNP}.
\end{itemize}
Note that while ``no'' instances in \textsc{coNP} admit easily verifiable certificates, this does not imply that these problems are easy to \emph{solve}; the intrinsic difficulty lies in discovering such a certificate or proving its absence. This is why both \textsc{coNP}-complete and \textsc{coNP}-hard problems are widely believed to be intractable.

\subsection{Polytope Representations}
A polytope $P \subseteq \mathbb{R}^n$ can be represented in two ways~\cite{Ziegler1995}:

\begin{enumerate}
\item The halfspace representation, also known as \textit{H-representation}, is given by
\[
P \ = \ \left\{ \ x \in \mathbb{R}^n \ \middle| \ H x \leq h \ \right\},
\]
where $H \in \mathbb{R}^{\mathsf{f} \times n}$ denotes a matrix and $h \in \mathbb{R}^\mathsf{f}$ a vector. In general, $P$ is a polyhedron, but if we additionally require that $H x \leq 0$ implies $x =0$, the set $P$ is bounded and, hence, a polytope.

\item The vertex representation, also known as \textit{V-representation}, is given by
\[
P \ = \ \operatorname{conv}\{ v_1, \dots, v_\mathrm{v} \} \ = \ \operatorname{conv}(V),
\]
the convex hull of a finite set of given points $v_1,\ldots,v_\mathrm{v} \in \mathbb R^n$. These points are often stored as columns of the matrix $V = [v_1,\ldots,v_\mathrm{v}]$.
\end{enumerate}
The polytope containment problem asks whether one polytope is contained in another. Of particular interest is the H-in-V containment problem: given an H-polytope $P$ and a V-polytope $Q$, with rational data,
\[
(H,h,V) \in \mathbb Q^{\mathsf{f} \times n} \times \mathbb Q^\mathsf{f} \times \mathbb Q^{n \times \mathsf{v}},
\]
decide whether $P \subseteq Q$.

\begin{theorem}[Freund and Orlin, 1985]
\label{thm::HinV}
The H-in-V polytope containment problem is \textsc{coNP}-complete.
This hardness remains even when $P$ is restricted to be a standard cube.
\end{theorem}

Note that the statement of the above theorem is an immediate consequence of~\cite[Thm.~1]{Freund1985}. The case that $P$ is restricted to a standard cube can be found in~\cite{Kellner2016}.

\subsection{Control Invariance}
Let us consider the discrete-time linear system
\begin{align}
\label{eq::System}
x^+ = A x + B u,
\end{align}
with state $x \in \mathbb{R}^n$, control $u \in \mathbb{R}^m$, and polytopic control constraint set $\mathbb U \subseteq \mathbb{R}^m$.

\begin{definition}
A set $X \subseteq \mathbb{R}^n$ is \emph{control invariant} for System~\eqref{eq::System} if
\[
\forall x \in X, \quad \exists u \in \mathbb U: \ A x + B u \in X.
\]
\end{definition}

Throughout this paper, we introduce the following standard convention~\cite{Schrijver1986}.

\begin{assumption}[Input Data]
\label{ass::InputData}
Following the standard convention in computational complexity, we assume that all input data—the matrices $A$ and $B$, the halfspace representation $(H,h)$ of the candidate polytope $X$, and the polytopic input constraint $\mathbb U$ (in either H- or V-representation)---have rational entries. This is a technical requirement that allows the formal definition of input size and polynomial-time reductions. The geometric results, however, naturally extend to real-valued data.
\end{assumption}

Given the above convention, the classical result by Gutman and Cwikel~\cite{Gutman1986} states that deciding whether a V-polytope is control invariant can be done in polynomial time. Namely, the V-polytope $P = \operatorname{conv}(V)$ is control invariant if and only if there exists a matrix $U = [U_1,\ldots,U_{\rm v}] \in \mathbb U \times \ldots \times \mathbb U$ and a non-negative matrix $\Lambda \in \mathbb R_+^{\mathrm{v} \times \mathrm{v}}$ with
\[
A V + B U \ = \ V \Lambda \ \ \text{and} \ \  \Lambda^\tr \ones = \ones \ \ \text{with} \ \ \ones = (1,\ldots,1)^\tr.
\]
Hence, as long as the control constraint set $\mathbb U$ is a given polytope and all data is rational, the above condition can be checked in polynomial time by using an LP solver~\cite{Khachiyan1979}; see also~\cite{Blanchini2015,Houska2025}.

\section{Verification of Control Invariance}

This section establishes the main complexity-theoretic result of the paper, which relies on the input data convention from Assumption~\ref{ass::InputData}.

\begin{theorem}
\label{thm::main}
The problem of verifying whether an H-polytope is control invariant for a discrete-time linear system with polytopic input constraints is \textsc{coNP}-complete. The hardness result holds regardless of whether the input constraint set $\mathbb U$ is specified in H- or V-representation.
\end{theorem}

The proof proceeds in two steps. We first establish \textsc{coNP}-hardness by a polynomial-time reduction from the H-in-V polytope containment problem. We then show membership in \textsc{coNP} by providing a succinct certificate for ``no'' instances.

\subsection{Hardness}
In order to construct a reduction from the H-in-V polytope containment problem, we set $P = [-1,1]^n$ to be the unit cube, and we set $Q = \operatorname{conv}(V)$ to be an arbitrary V-polytope. We consider System~\eqref{eq::System} with
\begin{align*}
A & \ = \ \left(
\begin{array}{cc}
0 & 0 \\
I & 0
\end{array}
\right), \ \ B \ = \ \left(
\begin{array}{cc}
0  \\
-V
\end{array}
\right), \\[0.16cm]
\text{and} \ \ \mathbb U & \ = \ \left\{ \ u \in \mathbb R_+^{m} \ \middle| \ \ones^\tr u = 1 \ \right\}.
\end{align*}
Now, we claim that the H-polytope $X = P \times \{ 0 \}$ is control invariant if and only if $P \subseteq Q$. In order to prove this, we start with any $p \in P$ and set $x = (p^\tr,0)^\tr \in X$. The successor is
\[
x^+ \ = \ \left(
\begin{array}{c}
0 \\
p - V u
\end{array}
\right).
\]
We claim that $X$ is control invariant if and only if $P \subseteq Q$. Now, we have $x^+ \in X$ if and only if $p - Vu = 0$, where $u$ must be an element of the unit simplex $\mathbb U$. But this means that $p = Vu$ is a convex combination of the columns of $V$. Hence,
\[
\exists u \in \mathbb U: \ p = Vu \iff p \in Q.
\]
In other words, for every $p \in P$, we must have $p \in Q$, which is exactly $P \subseteq Q$. This proves our claim that the H-polytope $X$ is control invariant if and only if $P \subseteq Q$.

Thus, a polynomial-time algorithm for verifying control invariance of H-polytopes would imply a polynomial-time algorithm for the H-in-V containment problem, which is \textsc{coNP}-hard. Moreover, the input constraint set $\mathbb U$ used in this reduction is a simplex, which has both a compact H-representation with $m+2$ inequalities and a compact V-representation with $m$ vertices. Hence, the reduction can be carried out in polynomial time regardless of whether the input constraint set is specified in H- or V-representation. This establishes \textsc{coNP}-hardness of the control invariance verification problem.

\begin{remark}
The above reduction may allow more inputs than states ($m > 2n$). This limitation is immaterial: one can append decoupled zero-state dimensions to $A$ and $B$ until the state dimension exceeds the input dimension, leaving the dynamics unchanged. Hence, the \textsc{coNP}-hardness result holds even under the restriction $m \le n$.
\end{remark}

\begin{remark}
The above reduction may produce a non-controllable system and a lower-dimensional candidate set. Neither restriction is essential. Standard perturbation arguments for rational polyhedra imply that the slice
\[
X = P \times \{0\}
\]
may be replaced by the full-dimensional polytope
\[
X_\varepsilon = P \times [-\varepsilon,\varepsilon]^n,
\]
where $\varepsilon>0$ is a sufficiently small rational number, without changing the outcome of the reduction~\cite{Schrijver1986}. Likewise, since controllability is an open and dense property of linear systems~\cite{Tchon1983}, the pair $(A,B)$ may be perturbed by an arbitrarily small rational amount to obtain a controllable system while preserving the reduction. Hence, the \textsc{coNP}-hardness result remains valid even under the additional assumptions of controllability and full-dimensionality.
\end{remark}

\subsection{Completeness}
It remains to show that the problem lies in \textsc{coNP}. For a ``no'' instance, a succinct certificate consists of a state $x \in X$ and a proof that no admissible input $u \in \mathbb U$ satisfies $A x + B u \in X$. 
This state \(x\) may be chosen to be a
vertex of \(X\), because the set of states in \(X\) admitting a feasible
successor is convex, so if all vertices of \(X\) admit such a
successor, then every point of \(X\) would.
Now suppose that $X$ is given in the form
\[
X \ = \ \{ x \in \mathbb R^n \mid H x \leq h \}.
\]
Then, the condition ``no $u \in \mathbb U$ satisfies $A x + B u \in X$'' is equivalent to the infeasibility of
\[
H(A x + B u) \le h, \quad u \in \mathbb U.
\]
By Farkas' lemma, infeasibility of this system has a dual certificate $\lambda$ (a vector of nonnegative multipliers) of size polynomial in the input data. But a deterministic verifier can check in polynomial time that
\begin{enumerate}
\item $x \in X$ by verifying $Hx \leq h$; and
\item the Farkas certificate $\lambda \geq 0$  proves that the system $H(Ax+Bu) \leq h,\; u \in \mathbb U$ has no solution~\cite{Schrijver1986}.
\end{enumerate}
This is because both checks consist of linear algebra operations (matrix-vector multiplications and comparisons) and are therefore polynomial-time computable.

Hence, every ``no'' instance admits a polynomial-size certificate that can be verified in polynomial time. Therefore, the control invariance verification problem is in \textsc{coNP}. Together with the hardness result from the previous subsection, this completes the proof of Theorem~\ref{thm::main}.

\begin{remark}
\label{rem::certificate}
The above coNP membership proof relies on the standard assumption that the certificate $(x,\lambda)$ have polynomial bit-length. This follows from well-known results on rational polyhedra: if the system $H(Ax+Bu) \le h,\ u \in \mathbb{U}$ is infeasible and $x$ is chosen as a rational vertex of $X$, then both $x$ and the Farkas certificate $\lambda$ have numerators and denominators bounded by determinants of submatrices of the input data~\cite{Schrijver1986}. Hence, the certificate size is polynomial in the input size.
\end{remark}

\section{Practical Implication}\label{sec::MCIpoly}

The complexity result of Theorem~\ref{thm::main} has immediate implications for the computation of maximal control invariant sets. In particular, the hardness of verifying control invariance transfers directly to any exact procedure for computing the maximal control invariant subset of a prescribed admissible polytope, under polynomial-time Turing reductions~\cite{Arora2009,Turing1939}.

\begin{cor}\label{cor::MCIset}
Suppose that the assumptions of Theorem~\ref{thm::main} hold, and let
$X\subseteq\mathbb{R}^n$ be an admissible polytope given by a rational
H-representation. The function problem of computing the maximal control
invariant subset of $X$ is \textsc{coNP}-hard under polynomial-time Turing 
reductions, for any exact output representation in which equality with a
rational H-polytope can be decided in polynomial time.
\end{cor}

\subsection{Proof of Corollary~\ref{cor::MCIset}}
Let $C(X)$ denote the maximal control invariant subset of
$X$. By definition, $X$ is control invariant if and only if
$C(X)=X$. If $X$ is control invariant, then its maximal control invariant subset
must be $X$ itself. Conversely, if $C(X)=X$, then $X$ is
control invariant because $C(X)$ is control invariant.

Suppose that an algorithm computes $C(X)$ exactly, for
every instance in the problem class considered in
Theorem~\ref{thm::main}, in time polynomial in the input size. If its output
representation permits equality with $X$ to be tested in polynomial time,
then control invariance of $X$ can be verified by first computing
$C(X)$ and then checking whether $C(X)=X$.
This gives a polynomial-time Turing reduction from the
\textsc{coNP}-hard control-invariance verification problem to the exact
computation of the maximal control invariant subset. Hence, the latter
function problem is \textsc{coNP}-hard under polynomial-time 
reductions.

\section{Conclusion}
We have proven that verifying control invariance of an H-polytope is coNP-complete. This result delineates a sharp complexity boundary between the tractable V-polytope case and the intractable H-polytope case. In model predictive control, where terminal sets are commonly specified in H-representation and constructed using a precomputed stabilizing control law, our result implies that any general verification procedure for such sets is inherently hard. Moreover, Corollary~\ref{cor::MCIset} implies that computing maximal 
control invariant terminal regions is \textsc{coNP}-hard. 
This establishes a fundamental limit on what can be achieved without additional structural assumptions.

\bibliographystyle{plain}
\bibliography{bib/references}

\end{document}